\documentclass[11pt]{elsarticle}
\usepackage[centertags]{amsmath}
\usepackage{amssymb}
\usepackage{amsthm}
\usepackage{amsfonts,amssymb}
\usepackage{enumerate}
\usepackage{mathrsfs}
\usepackage[nodots]{numcompress}
\usepackage{graphicx}

\journal{}

\theoremstyle{plain} 

\newtheorem{thm}{Theorem}

\theoremstyle{definition}

\theoremstyle{remark}
  
\newtheorem*{case1.1}{\textbf{\small Case 1.1}}

\DeclareMathOperator{\diverg}{Div}

\DeclareMathOperator{\pd}{\partial}

\DeclareFontFamily{OT1}{pzc}{}
\DeclareFontShape{OT1}{pzc}{m}{it}{<-> s * [1.10] pzcmi7t}{}
\DeclareMathAlphabet{\mathpzc}{OT1}{pzc}{m}{it}

\DeclareMathSymbol{\R}{\mathalpha}{AMSb}{"52}
\DeclareMathSymbol{\C}{\mathalpha}{AMSb}{"43}

\newcommand{\set}[1]{\left\{#1\right\}}

\newcommand{\comment}[1]{}

\newcommand{\mscs}{\mathscr{S}}

\newcommand{\bv}{\mathbf{v}}

\newcommand{\beq}{\begin{equation}}
\newcommand{\eeq}{\end{equation}}
\newcommand{\beqS}{\begin{equation*}}
\newcommand{\eeqS}{\end{equation*}}
\newcommand{\balign}{\begin{align}}
\newcommand{\ealign}{\end{align}}
\newcommand{\bsube}{\begin{subequations}}
\newcommand{\esube}{\end{subequations}}

\begin{document}

\begin{frontmatter}



\title{Lie group classification and conservation laws of a class of hyperbolic equations}

\author{J.C. Ndogmo\corref{cor1}}
\ead{jean-claude.ndogmo@univen.ac.za}

\address{Department of Mathematics and Applied Mathematics\\
University of Venda\\
P/B X5050, Thohoyandou 0950, South Africa}

\cortext[cor1]{Corresponding author}


\begin{abstract}
A new method for the Lie group classification of differential equations is proposed. It is based of the determination of all possible cases of linear dependence of certain indeterminate appearing in the determining equations of symmetries of the equation. The method is simple and systematic and applied to a family of hyperbolic equations. Moreover, as the said family  contains several known equations with important physical applications,  low-order conservation laws of some relevant equations from the family are computed, and the results obtained are discussed with regard to the symmetry integrability of a particular class from the underlying family of hyperbolic equations.
\end{abstract}


\end{frontmatter}

\section{Introduction}
\label{s:intro}
The Lie group classification of differential equations, which consists in determining all symmetry classes admitted by an equation according to the values of the parameters or arbitrary functions labelling the given family of differential equations, has been carried out in the literature only in a more or less ad hoc manner \cite{ovsyC, gungorC, moyoC2, moyoC3, ndogC, BihloC, ZdanovC, bagderC}. This has resulted as pointed out in \cite{ bagderC}  in a number of revised classification results published in the literature being wrong. Attempts  to find  a somewhat systematic method for this classification problem has however been made. This includes the well-known algebraic method which can be traced back to Lie's work on symmetry algebras of ordinary differential equations ({\sc ode}s), and which have been upgraded or applied to several papers \cite{34ngang, moyoC3, 2ngang, 13ngang, gungorC, gungorC2}. Another similar attempt for a systematic method  was proposed for linear {\sc pde}s in \cite{ZdanovC}. All these suggested method require in general a lot of computations and analysis even for simpler cases of equations. They are also limited by the availability of the large amount of information they usually require. For instance, the algebraic method can only be carried out when any symmetry algebra of the given family of equations falls, in general,  within the existing classification of low dimensional Lie algebras.\par

In this paper, we proposed a new systematics method for the Lie group classification of  differential equations. It is based on the determination of all possible cases of linear dependence of certain indeterminate appearing in the determining equations for the Lie point symmetry algebra of the given family of equations.  The method is presented through an application to the family of scalar hyperbolic equations
\begin{equation}\label{eq:M}
u_{xy} = F (u, u_x),   
\end{equation}
labelled by the arbitrary function $F=F(u, u_x),$ where $u=u(x,y),$ and a variable in a subscript denotes partial derivative w.r.t. the variable, so that $u_x= \pd u/ \pd x,\, u_y= \pd u/ \pd y, \, u_{xy}= \pd^2 u/ (\pd x\, \pd y)$ and so on.  This systematic  method turns out to be simple and yields a relatively fast classification of the family \eqref{eq:M} of equations. \par

Given that the family  \eqref{eq:M} contains several well-known equations with important physical applications, including amongst others Liouville's equation \cite{popLv, winterLv1, winterLv2} and the sine-Gordon equation, low-order conservation laws of a number of relevant equations from this family are computed using the direct method \cite{olv93, muatjetjeja, anco_kara, anco_bluman}. Indeed, the sine-Gordon equation occurs for instance in the study of surfaces of constant negative curvature as well as in the study of crystal dislocations, and it's solutions  possess some soliton properties \cite{sineG1, sineG2}. Liouville's equation on the other hand appears in the study of isothermal coordinates in classical mechanics and statistical physics \cite{dubrovin}. Some of these computations correspond to the case of arbitrary labelling functions. Although some of the conservation laws found are finite in number, they are infinite for other equations, and the connection between such equations and their symmetry integrability is discussed. It will be assumed that $F$ is a nonconstant function and either $F=F(u),$ or $F=F(u_x).$

\section{Lie group classification} \label{s:classif}
   The Lie group classification of \eqref{eq:M} will be performed under point transformations.  We denote by $\bv= \xi \pd_x + \eta \pd_y + \phi \pd_u,$ a generic symmetry vector of \eqref{eq:M}, where $\xi, \eta$ and $\phi$ are functions of $x,y$ and $u.$ By generic symmetry vector of \eqref{eq:M} we mean a linear combination of all symmetry generators of \eqref{eq:M}. \par

   The determining equations (D) for the symmetries of \eqref{eq:M}, whose expression is omitted here, show that the following conditions should hold.

\begin{subequations}\label{eq:det1}
\begin{align}
\xi_u&=\xi_y= \eta_u= \eta_x= \phi_{uu}=0,\\
\intertext{Consequently,}
\xi&=\xi(x),\quad \eta= \eta(y),\quad \text{ and } \phi= g\, u+ h
\end{align}
\end{subequations}
 for some functions $g=g(x,y)$ and $h= h(x,y).$  Given the conditions on the function $F$ in \eqref{eq:M}, we will treat separately two cases.
\subsection{Case 1: F=F(u).}

Updating (D) with this expression for $F$ and the new values for $\xi, \eta$ and $\phi$ given by \eqref{eq:det1} shows that $g_x=g_y=0$ must hold. Hence $g= A$ must be a constant. With this value for $g$ the other remaining equation in (D) reduces to
\begin{equation} \label{eq:det2}
- F_u h - u F_u A  + F(A- \eta_y - \xi_x) + h_{xy}=0.
\end{equation}
Before proceeding any further with the ongoing analysis, it is appropriate at this point to determine the equivalence group of \eqref{eq:M} corresponding to the current expression for $F.$
\begin{thm}\label{th:eqvF(u)}
The group of equivalence transformations of the hyperbolic equation
\begin{subequations}\label{eq:Fu}
\begin{align}
  u_{xy} &= F(u) \label{eq:Fu1} \\[-3mm]
\intertext{ consists of linear transformations} \vspace{-4mm}
  x  &= a t+ b,\quad y= c z+ d,\quad  u= r w + s,  \label{eq:Fu2} \\[-3mm]
\intertext{for some arbitrary constants $a, b, c, d, r, s,$ where $w=w(t,z).$ The corresponding transformed equation has expression}\vspace{-4mm}
  w_{tz} &=  H(w),  \label{eq:Fu3}
\end{align}
where $H= \frac{ac}{r} F(r w+ s)$.
\end{subequations}
\end{thm}

\begin{proof}
  Subjecting \eqref{eq:Fu1} to the most general invertible point  transformation
\begin{equation}\label{eq:gnlpoint}
 x= \sigma(t,z,w),\quad y= \tau(t,z, w), \quad u=\rho(t,z,w),
\end{equation}

and requesting that the resulting transformed equation be of the form \eqref{eq:Fu3} for some arbitrary function $H$ leads to the vanishing of all other terms in the transformed equation.  The resulting constraints on the functions $\sigma, \tau,$ and $\rho$ reduce the admissible transformations of $x,y$ and $u$ to \eqref{eq:Fu2}, and the transformed equation to \eqref{eq:Fu3} with the specified expression for $H.$ This completes the proof.
\end{proof}
As a consequence of Theorem \ref{th:eqvF(u)}, in the determination of the Lie group classification of \eqref{eq:Fu1}, one can always replace, in particular, the function $F=F(u)$ by any of its nonzero scalar multiple, although admissible transformations of $F$ are quite more general than this. Moreover, for two given equations \eqref{eq:Fu1} and \eqref{eq:Fu3} equivalent under the transformation \eqref{eq:Fu2}, we may write $F\sim H,$ and it is clear that the relation $\sim$ thus defined is an equivalence relation on the space of functions $F= F(u).$  Keeping this in mind, we will now proceed to solving Equation \eqref{eq:det2}.\par

Our general procedure for the Lie group classification problem goes as follows. In each term of \eqref{eq:det2}, we consider  those factors involving only functions of $u$ as indeterminate in a polynomial expression. In that way the vanishing coefficients of the polynomial expression can be unambiguously determined if we know exactly which of the indeterminate are nontrivially linearly dependent and which ones are not. Here, a nontrivial linear dependence is one for which all coefficients in the corresponding vanishing linear combination are nonzero. This will be achieved by assuming that there are exactly $m$ linearly dependent indeterminate, with $ 0\leq m\leq n,$ where $n$ is the total number of indeterminate in the equation. For each such value of $m$ all possible subsets of $m$ linearly dependent indeterminate will be considered, which amounts to requesting that their linear combination with some arbitrary nonzero  constant coefficients vanishes. Alternatively,  the Wronskian of the $m$ indeterminate should vanish, under the assumption that $F$ is a smooth function. Each set of $m$ indeterminate will thus yield a differential equation to be satisfied by $F$ and hence a corresponding value for $F.$ For each value of $m$ there should be at most $\binom{n}{ m}$ valid differential equations, where $\binom{n}{ m}$ is the usual binomial coefficient. Some of the values of $F$ thus achieved for all possible values of $m$ should also be compared to single out cases where redundant linearly dependent indeterminate occur.  Such comparisons will usually be done by mere inspection.  Symmetries determined by each function $F$  in this updated list will then give rise to a complete list of symmetry classes for the equation.\par

For the study of Equation \eqref{eq:det2}, we first consider the case where the function $F$ is arbitrary, which according to the above notation corresponds to having $m=0$. The vanishing of the coefficient of $F_u$ then shows that $h=A=0,$ while the vanishing of the coefficient of $F$ yields $\xi_x= - \eta_y.$ Consequently, for $F$ arbitrary the symmetry algebra of \eqref{eq:Fu1} is generated by
\begin{equation}\label{eq:FuArb}
 \bv = (k_1 x + k_2) \pd_x + (k3-k_1 y) \pd_y,
\end{equation}
where here and in the sequel the $k_j,\, \text{ for }j=1,2,\dots$ represent arbitrary constants. \par

Assuming now that $F$ is not arbitrary in  \eqref{eq:det2}, the set  $\mscs=\set{F_u, uF_u, F, 1}$ of indeterminate consists  precisely of $n=4$ elements, and since none of these vanishes one should assume that $m \geq 2.$ \par

\subsubsection{Two indeterminate are linearly dependent.}
The constraints on $F$ determined by the linear dependence of all possible subsets of $m=2$ elements from $\mscs$  are given by
\begin{subequations}\label{eq:Fu:2}
  \begin{align}
    0= & F_u^2    \label{eq:Fu:2E1}\\
    0= & F_u^2 - F F_{uu}   \label{eq:Fu:2E2} \\
    0= & -F_{uu}    \label{eq:Fu:2E3} \\
    0= & u F_u^2  - F(F_u + u F_{uu})    \label{eq:Fu:2E4} \\
    0= & -F_u - uF_{uu}    \label{eq:Fu:2E5} \\
    0= & F_u    \label{eq:Fu:2E6}
  \end{align}
\end{subequations}
Since by assumption the function $F$ may not assume a constant value, it follows by inspection that \eqref{eq:Fu:2E1} and \eqref{eq:Fu:2E6} yield invalid solutions. Other values of $F$ yielding only the symmetry generator \eqref{eq:FuArb}  are also to be excluded from the updated list of solutions of \eqref{eq:Fu:2}. The relevant values of $F$ determined by \eqref{eq:Fu:2} and the corresponding generic symmetry generators $\bv$  are therefore given as follows.  Here and in the sequel the $a_j, \text{ for } j= 1,2,\dots$ are arbitrary constants. Hence we have,
\begin{subequations}\label{eq:Fu:2Sym}
  \begin{alignat}{2}
  \text{ For }F &\sim e^{ u},\quad&   \bv &= \xi(x) \pd_x + \eta(y) \pd_y -\left[ \xi_x(x) + \eta_y(y)\right]\pd_u.     \label{eq:Fu:2s1}\\
  \text{ For }F &\sim u,\quad&   \bv &= (k_1 x + k_2) \pd_x + (-k_1 y + k_3)\pd_y + (k_4\, u + h)\pd_u ,       \label{eq:Fu:2s2} \\
  \intertext{where $h$ is a solution of the original equation \eqref{eq:Fu1}, while $\xi$ and $\eta$ are arbitrary functions of their arguments.}
  \text{ For }F &\sim {u^{a_1}},&               &  a_1 \neq 1,  \notag\\[-3mm]
                &               &           \bv &= (k_1 x + k_2)\pd_x + (k_3 y + k_4)\pd_y - \frac{(k_1 + k_3)u}{a_1 -1}\pd_u.    \label{eq:Fu:2s3}
  \end{alignat}
\end{subequations}

\subsubsection{Three indeterminate are linearly dependent.}
We now move on to consider the case where three elements in $\mscs$ are linearly dependent. The conditions on $F$ imposed by all possible such sets of three linearly dependent indeterminate are given by the following equations.
\begin{subequations}\label{eq:Fu:3}
  \begin{align}
    0= & -F_u^2 F_{uu} + 2 F F_{uu}^2 - F F_u F_{uuu}  \label{eq:Fu:3E1}\\
    0= & 2F_{uu}^2 - F_u F_{uuu}   \label{eq:Fu:3E2} \\
    0= & F_{uu}^2 - F_u F_{uuu}   \label{eq:Fu:3E3} \\
    0= &  u F_{uu}^2 - F_u (F_{uu}+ u F_{uuu})    \label{eq:Fu:3E4}
  \end{align}
\end{subequations}

The change of variables $F= e^{w}$ transforms \eqref{eq:Fu:3E1} into
\begin{equation}\label{eq:Fu:3E1A}
2 w_{uu}^2 - w_u w_{uuu}=0,
\end{equation}
and the latter  has general solution $w= b_2 \ln(u+ b_1) + b_3,$ where here and in the sequel the $b_j,\, \text{ for $j= 1,2,\dots$}$ denote arbitrary  constants. Consequently the general solution of \eqref{eq:Fu:3E1} is given by
\begin{subequations}\label{eq:Fu3E1Sol}
\begin{align}
F   &= a_3 (u+ a_1)^{a_2} \sim (u+ a_1)^{a_2} \label{eq:Fu3E1SolF}\\
\intertext{with corresponding symmetry}
\bv &= (k_1 x + k_2)\pd_x + [k_3 - (k_1 + k_4 (a_2-1))y ] \pd_y + k_4 (a_1 + u)\pd_u.  \label{eq:Fu3E1SolV}
\end{align}
\end{subequations}

The other equations in \eqref{eq:Fu:3} are solved in a way more or less similar to that for Equation \eqref{eq:Fu:3E1}. However they yield solutions such as
\[
F= a_2 \ln (u+ a_1) + a_3,\quad \text{ or } F= \frac{a_2e^{a_1 u}}{ a_1} + a_3,\quad \text{ or } F= \frac{a_2 u^ {a_1}}{ a+1} + a_3
\]
which are often more general solutions then those found for $m=2,$  but which in any case do not yield symmetries other than those given by \eqref{eq:FuArb}. \par

\subsubsection{Four indeterminate are linearly dependent.}
We now consider the case where there are four linearly dependent elements in $\mscs.$ There is only one such possibility consisting of the whole set $\mscs$ and the corresponding  constraint on $F$ is given by the equation
\begin{equation}\label{eq:GivFu2}
  F_{uu}^2 F_{uuu} - 2 F_u F_{uuu}^2 + F_u F_{uu} F_{uuuu}=0.
\end{equation}
Setting $F= \int G(u) du$ and then $G= e^{w(u)}$ reduces \eqref{eq:GivFu2} to
\[
-2w_{uu}^2 + w_u w_{uuu}=0,
\]
so that $w= c_3 + c_2 \ln(u + c_1),$  $G= b_3 (u+ b_1)^{b_2},$ for some arbitrary constants $ c_j, b_j, \text{ where } j=1,\dots, 3.$ Consequently, for $b_2 \neq -1,$ $F \sim  (u+ a_1)^{a_2} + a_3.$  We note that $F$ is linear in this case if and only if $a_2=1,$ in which case $F\sim u$ corresponds to the already solved case \eqref{eq:Fu:2s2}.
On the other hand, for  $a_2 \neq 1$ and  $a_3=0,$ $F$ corresponds to the solved case \eqref{eq:Fu3E1Sol}, while for $a_2 \neq 1,$ and  $a_3\neq 0,$ $F$  yields the  symmetry  \eqref{eq:FuArb} corresponding to arbitrary functions.

Similarly, for  $b_2= -1$ in the expression of $G= b_3 (u+ b_1)^{b_2},$ one gets $F \sim  \ln (u+ a_1) + a_2$ which as already seen also yields a symmetry generator given by \eqref{eq:FuArb}.

\begin{thm}\label{th:classFu}
Denote by $L$ the symmetry algebra of the hyperbolic
 equation $u_{xy}= F(u)$ and by $\bv$ the generic symmetry vector in $L.$
\begin{enumerate}[\rm{(}a\rm{)}]
\item For $F(u) \sim u,$ $L$ is infinite dimensional and
\[
\bv = (k_1 x + k_2) \pd_x + (-k_1 y + k_3)\pd_y + (k_4\, u + h)\pd_u
\]
where $h$ is a solution of the original equation \eqref{eq:Fu1}.
\item For $F(u) \sim   e^{ u} $, $L$ is infinite dimensional and
\[
\bv =  \xi(x) \pd_x + \eta(y) \pd_y -\left[  \xi_x(x) + \eta_y(y) \right] \pd_u.
\]
\item  For $F(u) \sim (u+ a_1)^{a_2},$ $a_2\neq 1,$ $L$ has dimension $4$ and
\[
\bv = (k_1 x + k_2)\pd_x + [k_3 - (k_1 + k_4 (a_2-1))y ] \pd_y + k_4 (a_1 + u)\pd_u.
\]
\item  For $F(u) \sim u^{a_1},$ with $a_1\neq 1,$ $L$ has dimension $4$ and
\[
\bv =  (k_1 x + k_2)\pd_x + (k_3 y + k_4)\pd_y - \frac{(k_1 + k_3)u}{a_1 -1}   \pd_u
\]
\item For any other function $F(u),$ $L$ has dimension $3$ and
\[
 \bv = (k_1 x + k_2) \pd_x + (k3-k_1 y) \pd_y.
\]
The five symmetry classes thus obtained are pairwise nonequivalent  and make up all possible symmetry classes of \eqref{eq:Fu1}.
\end{enumerate}
\end{thm}

\begin{proof}
The statements (a)-(e) as well as the fact that the listed symmetry classes exhaust all possible symmetry classes of \eqref{eq:Fu1} are just a summary of the results established immediately before the theorem, and it only remains to show that the stated symmetry classes are nonequivalent. The fact that the symmetry algebra in Case (e) has dimension $3$ which is less than the dimension of $L$ in all other cases establishes the non equivalence  between the symmetry algebra in Case (e) with those in all other cases. For the cases (a)-(d), the corresponding symmetry algebras are nonequivalent, precisely because in view of the equivalence transformations \eqref{eq:Fu2}, the associated functions $F,$ and hence the associated equations of the form \eqref{eq:Fu}, are nonequivalent. This completes the proof of the Theorem.
\end{proof}

\comment{
It  should be noted that the solution \eqref{eq:A=0e3} is a limiting case of the solution $a_4(u+ a_1)^{a_2} + a_3$ obtained for \eqref{eq:GivFu2}, and in fact it solves\eqref{eq:GivFu2}. Indeed, we may write
\[
(a_1+ u)^{a_2}=??
\]
}
\subsection{Case 2: $F=F(u_x).$}\label{subsecFux}
It should be noted that using a simple symmetry argument based on the structure of the resulting equation $u_{xy}= F(u_x)$ with respect to the variables $x,$ and  $y,$ the results from this section implicitly include those corresponding to the case $F=F(u_y.)$ We begin the Lie group classification in the actual case by finding the corresponding Lie group of equivalence transformations.
\begin{thm}\label{th:eqvFux}
The group of equivalence transformations of the class of hyperbolic equations
\begin{subequations}\label{eq:Fux}
\begin{align}
  u_{xy} &= F(u_x) \label{eq:Fux1} \\[-3mm]
\intertext{ consists of linear transformations} \vspace{-4mm}
  x  &= a t+ b,\quad y= c z+ d,\quad  u= r w + S(z),  \label{eq:Fux2} \\[-3mm]
\intertext{for some arbitrary constants $a, b, c, d, r,$ and arbitrary function $S=S(z),$ where $w=w(t,z).$ The corresponding transformed equation has expression}\vspace{-4mm}
  w_{tz} &=  H(w_t),  \label{eq:Fux3}
\end{align}
where $H= \frac{a\,c}{r} F(\frac{r\, w_t}{a})$.
\end{subequations}
\end{thm}

\begin{proof}
Under the general point transformation \eqref{eq:gnlpoint}, the vanishing of the coefficient of $w_{tt}$ and $w_{zz}$ in the transformed version of \eqref{eq:Fux1} shows that either $\sigma= \sigma(t) \text{ and } \tau= \tau(z)$ or  $\sigma= \sigma(z) \text{ and } \tau= \tau(t)$ must hold. But due to the fact that the arbitrary function $F$ should not depend explicitly on $z,$  only the first possibility for $\sigma$ and $\tau$ may hold. Therefore, applying the transformation \eqref{eq:gnlpoint} with $\sigma= \sigma(t) \text{ and } \tau= \tau(z)$ and requesting that the transformed equation keeps the form \eqref{eq:Fux1}, and in particular that the argument of the function $F$ in the transformed equation depends only on $w_t$ immediately yields the stated results. This completes the proof of the Theorem.
\end{proof}

It should be noted that the transformations \eqref{eq:Fux2} are linear and with almost only constant coefficients and are therefore quite weak. In particular, under these transformations, the functions $F=u$ and $F= u+ c$ are not equivalent for a nonzero value of the constant $c.$

When the expression of $F$ in \eqref{eq:M} is reduced to $F=F(u_x),$ the corresponding determining equations (D) show that in addition to the conditions \eqref{eq:det1} restricting the components $\xi, \eta,$ and $\phi$ of the symmetry generators one should also have $g=g(y)$ in \eqref{eq:det1}. The resulting expression of the determining equations (D) is reduced to the single equation
\begin{equation}\label{eq:detFux}
u_x\, g_y + F(g - \eta_y - \xi_x) + u_x F_{u_x}(\xi_x - g)  - F_{u_x} h_x + h_{xy} =0
\end{equation}
which will be used to classify \eqref{eq:Fux1}.\par

When the function $F$ in \eqref{eq:detFux} is arbitrary and correspondingly $m=0,$ an analysis similar to the one done for the preceding  case $F=F(u)$ shows that the generic symmetry is given by
\begin{equation}\label{eq:FuxArb:sym}
\bv=  (k_1 x + k_2)\pd_x + k_3\pd_y + ( k_1 u + P)\pd_u
\end{equation}
where $P= P(y)$ is an arbitrary function. It is a remarkable fact that the class \eqref{eq:Fux1}   of equation has an infinite dimensional principal symmetry algebra. It might therefore be possible that this class of equations is linearizable by certain types of transformations, although it clearly follows from Theorem \ref{th:eqvFux} that \eqref{eq:Fux1}  is not linearizable by point transformations. \par

For the rest of the classification of \eqref{eq:Fux1}, we continue with the application of our method based on the determination of all possible cases of linearly dependent subsets of $m$  indeterminate in \eqref{eq:detFux}, with $2 \leq m \leq n,$   where in this case $n=5$ is the cardinality of the set $\mscs=\set{u_x, F, u_x F_{u_x}, F_{u_x},1}$ of indeterminate.    In the set of equations representing the constraints on $F$ for a particular value of $m,$ only those equations corresponding to  new and non redundant solutions will usually be represented.

\subsubsection{Two indeterminate are linearly dependent.}
Therefore, for $m=2$ the constraints on $F$ are given by
\begin{subequations}\label{eq:Fux:m2E}
\begin{alignat}{2}
   0 &= -F + {u_x}F_{u_x},                         &\qquad \qquad              0 &= - F_{u_x} + {u_x} F_{{u_x}{u_x}}       \\
   0 &={u_x}^2 F_{{u_x}{u_x}},                          &\qquad \qquad              0 &= -  F_{u_x}^2 + F F_{{u_x}{u_x}}    \\
 0 &= -{u_x}F_{u_x}^2 + F(F_{u_x} + {u_x}F_{{u_x}{u_x}}),         &\qquad \qquad               0 &= - F_{u_x} - {u_x} F_{{u_x}{u_x}}.
\end{alignat}
\end{subequations}
To simplify notation, we will make use of the following vector fields, where here and in the sequel, unless otherwise stated $\xi$ and $H$ are arbitrary functions of $x$  while  $g, S,$ and $P$ are arbitrary functions of $y,$ and $\alpha$ is an arbitrary constant.

 More precisely, for each value of $m$ we denote by $V_{m, j}$ the  $j$-th generic generator, as they consecutively occur,  of the symmetry class associated with a solution of an {\sc ode} representing the condition of linear dependence of $m$ indeterminate. Therefore, let us set
\begin{subequations}\label{eq:Fuxm2sym}
\begin{align}
 V_{2, 1} &= \xi \pd_x + (k_1+ g)\pd_y +  \left[g u + a_1 (x g-\xi)+ e^y H(x) + S(y)\right]\pd_u \label{eq:Fuxm2sym1}\\
 \begin{split}\label{eq:Fuxm2sym2}
 V_{2, 2} &=   (k_5+k_6 x)\pd_x +   (k_2+y (k_3+k_4 y))\pd_y   \\
    &\quad + [k_1+P-k_4 x+u (-k_3+k_6-2 k_4 y) ]\pd_u
 \end{split}\\
 \begin{split}\label{eq:Fuxm2sym3}
 V_{2, 3} &= (k_1+k_2 x)\pd_x + \frac{2 \sqrt{\alpha} k_5-k_4\cos(2 \sqrt{\alpha} y)+k_3\sin(2 \sqrt{\alpha} y)}{2
\sqrt{\alpha}}\pd_y          \\
      &+ \left[ P+k_2 u-\left(k_3 u+\sqrt{\alpha} k_4 x\right)\cos(2 \sqrt{\alpha} y)+\left(-k_4 u+\sqrt{\alpha} k_3 x\right)
\sin(2 \sqrt{\alpha} y)\right]\pd_u
 \end{split}\\
  V_{2, 4} &= (k_3+k_4 x)\pd_x   + (k_1+k_2 y)\pd_y + \left[P+u \left(k_4+\frac{k_2}{1-a_1}\right) \right] \pd_u \label{eq:Fuxm2sym4}\\
  V_{2, 5} &= (k_2+k_3 x)\pd_x  +   (k_1-k_5 y)\pd_y  +  (k_4+P+k_3 u+k_5 x)\pd_u. \label{eq:Fuxm2sym5}
\end{align}
\end{subequations}
The solutions $F$ of \eqref{eq:Fux:m2E} represented by their canonical forms under \eqref{eq:Fux2} and the corresponding generic symmetry vector $\bv$ are given as follows.
\begin{subequations}\label{eq:Fux:m2Sym}
\begin{alignat}{2}
 F & \sim {u_x} + a_1,                           &\qquad             \bv &=  V_{2, 1}      \label{eq:Fux:m2Sym1} \\
 F & \sim {u_x}^2,                               &\qquad             \bv &=  V_{2, 2}      \label{eq:Fux:m2Sym2} \\
 F &\sim   {u_x}^2+ \alpha, (\alpha \neq 0),             &\qquad             \bv &=  V_{2, 3}   \label{eq:Fux:m2Sym3} \\
 F  & \sim {u_x}^{a_1}, (a_1 \notin \set{1,2}),     &\qquad             \bv &=  V_{2, 4}      \label{eq:Fux:m2Sym4}\\
 F & \sim e^{u_x},                             &\qquad             \bv &=  V_{2, 5} \label{eq:Fux:m2Sym5}
\end{alignat}
\end{subequations}

\subsubsection{Three indeterminate are linearly dependent.}
In order to write down more concisely the set of conditions on $F$ corresponding to $m=3,$ we set $u_x=\nu.$ The equations are then given by
\begin{subequations}\label{eq:Fux:m3}
\begin{align}
 0 &= F_{\nu \nu \nu } \label{eq:Fux:m3E1}\\
 0 &=- F^2_{\nu \nu } + F_\nu  F_{\nu \nu \nu }   \label{eq:Fux:m3E2}\\
 0 &= -F_{\nu }^2 F_{\nu \nu }+2 F F_{\nu \nu }^2-F F_{\nu } F_{\nu \nu \nu } \label{eq:Fux:m3E3} \\
 0 &=-\nu  F_{\nu \nu }^2-F F_{\nu \nu \nu }+F_{\nu } \left(F_{\nu \nu }+\nu  F_{\nu \nu \nu }\right) \label{eq:Fux:m3E4}\\
 0 &= -F \left(2 F_{\nu \nu }+\nu  F_{\nu \nu \nu }\right)+\nu  \left(-\nu  F_{\nu \nu }^2+F_{\nu } \left(2 F_{\nu \nu }+\nu  F_{\nu \nu \nu }\right)\right). \label{eq:Fux:m3E5}
\end{align}
\end{subequations}
Here again we set
\begin{align}
\begin{split}
 V_{3,1} &=  (k_1+k_2 x)\pd_x + [k_3+y (k_4+k_5 y)] \pd_y \\
           &+ [P+u (k_2-k_4-2 k_5 y)-\frac{1}{2} x (\beta k_4 + 2 k_5+2 \beta k_5 y)]\pd_u
 \end{split}\\
\begin{split}
 V_{3,2} &=  (k_1+k_2 x)\pd_x +  \left[\frac{2 \delta k_5-k_4 \cos(2 \delta y)+k_3 \sin(2 \delta y)}{2 \delta}\right]\pd_y \\
    &+ \bigg[ P+k_2 u-\frac{1}{2} (2 k_3
u+\beta k_3 x+2 \delta k_4 x) \cos(2 \delta y)\\
     &-\frac{1}{2} (2 k_4 u-2 \delta k_3 x+\beta k_4 x) \sin(2 \delta y)\bigg]\pd_u
\end{split}\\
\begin{split}
 V_{3,3} &= (k_1+k_2 x)\pd_x +\left[ k_5+\frac{e^{2 \delta y} k_3-e^{-2 \delta y} k_4 }{ 2 \delta} \right]\pd_y\\
           & + \big[P+k_2 u+\frac{1}{2} e^{-2 \delta y} k_4 (-2 u+2 \delta x-\beta x)\\
& -\frac{1}{2} e^{2 \delta y} k_3 (2 u+(2 \delta+\beta) x)\big]\pd_u.
\end{split}
\end{align}
In view of the equivalence transformations \eqref{eq:Fu2}, the solution $F$ of \eqref{eq:Fux:m3E1} satisfies $F \sim u_x^2 + \beta u_x + \alpha.$ Denote by $\Delta = (\beta^2 - 4 \alpha )$  the discriminant associated with this quadratic polynomial and assumed to be  a real number, and set $\delta^2 = |\Delta|/4$ if $\Delta \neq 0,$ so that $\delta^2 >0.$  Then according to the sign of $\Delta$ the solution $F$ splits into three cases with corresponding generic symmetry vectors $\bv$ as follows.
\begin{subequations}\label{eq:Fux:m3Sym1}
\begin{alignat}{2}
   F &\sim   \beta^2/4         + \beta u_x + u_x^2,             &\qquad             \bv &=  V_{3,1}  \\
  F &\sim   (\beta^2/4 + \delta^2) + \beta u_x + u_x^2,             &\qquad             \bv &=  V_{3,2}   \\
   F  &\sim (\beta^2/4 - \delta^2) + \beta u_x + u_x^2,             &\qquad             \bv &=  V_{3,3}.
\end{alignat}
\end{subequations}
The general solution $F $ of \eqref{eq:Fux:m3E2} satisfies $F \sim e^{\alpha u_x} + \beta$ and for $\beta =0$ the corresponding symmetry is given by \eqref{eq:Fux:m2Sym5}, while for $\beta\neq 0,$ $F \sim e^{\alpha u_x} + 1$ and the corresponding symmetry is given by
\begin{equation} \label{eq:Fux:m3E2Sym}
 \bv  =  (k_1 + k_2 x) \pd_x + \left( k_4 - \frac{k_3 e^{- \alpha y}}{\alpha} \right)\pd_y + \left(P(y) + k_2 u - \frac{k_3 x e^{- \alpha y}}{\alpha} \right)\pd_u.
\end{equation}

The next equation to consider from \eqref{eq:Fux:m3} is \eqref{eq:Fux:m3E3} which is up to a renaming of the unknown function's argument the same as  \eqref{eq:Fu:3E1} already solved and whose solution satisfies $F\sim (\alpha + u_x)^\beta.$ The symmetries corresponding to the solution $F$ are given as follows for $\alpha \neq 0,$ the case $\alpha =0$ being already solved in \eqref{eq:Fux:m2Sym4}.
\begin{subequations}\label{eq:Fux:m3Sym3}
\begin{alignat}{2}
   F &\sim   (\alpha + u_x)^2,    (\alpha  \neq 0),         &\qquad             \bv &=  V_{3,4} \label{eq:Fux:m3Sym3a}\\
  F &\sim    (\alpha + u_x)^\beta, (\beta \notin \set{1,2}, \alpha \neq 0),             &\qquad             \bv &=  V_{3,5}, \label{eq:Fux:m3Sym3b}
\end{alignat}
\end{subequations}
where
\begin{subequations}
\begin{align}
\begin{split}
 V_{3,4} &=  (k_1 + k_2\, x)\pd_x + [ k_3 + y(k_4+ k_5\, y)] \pd_y \\
           &\quad+ [P(y) -x\, (\alpha\,  k_4 + k_5 + 2 \alpha\, k_5\, y)  + (k_2 - k_4- 2 k_5\, y)\,u]\pd_u
 \end{split}\\
\begin{split}
 V_{3,5} &=  (k_2 + k_3\, x)\pd_x +  \left[  k_1 - (\beta -1) \frac{k_5\, y}{\alpha} \right]\pd_y \\
    &\quad+ \big[ k_4 + k_5\, x + P(y) + \left(k_3 + \frac{k_5}{\alpha}\right)u   \big]\pd_u.
\end{split}
\end{align}
\end{subequations}

    Up to this point the constraints on $F$ resulting from the condition of linear dependence of indeterminate have been expressed typically as $m$-th order {\sc ode}s, and we have luckily been able to solve all such {\sc ode}s. However, some of them, such as \eqref{eq:Fux:m3E4}   and \eqref{eq:Fux:m3E5},  are hard to solve and we have to resort to the more direct method of expressing the condition of linear dependence as the vanishing of a linear combination of the indeterminate with some arbitrary constant coefficient. In that way all such constraints are mere  linear first order {\sc ode}s. The drawback with this method lies in the rigidity of arbitrary constants in the resulting general solution which albeit straightforward to find, are somewhat much harder to reduce to a simpler and suitable form, due to the fact that the coefficients in the linear combination are in fact  arbitrary parameters of the required solution.\par

    For Equation \eqref{eq:Fux:m3E4}, the corresponding vanishing of the linear combination of indeterminate takes the form
\begin{subequations}\label{eq:FuxLC1}
\begin{align}
0 &= \alpha u_x + \beta F +\sigma F_{u_x},\label{eq:FuxLC1E1}\\
\intertext{with solution}
F &= \delta e^{-\frac{\beta u_x}{ \sigma}} - \alpha(- \sigma + \beta u_x)/ \beta^2,  \label{eq:FuxLC1E2}
\end{align}
\end{subequations}
where the constant of integration $\delta$ together with the coefficient $\alpha, \beta,$ and $\tau$ are to be considered as arbitrary parameters in the solution  \eqref{eq:FuxLC1E2}.  With a bit of manipulation, the latter solution can be put into the form
\[
F= \lambda (1+ \beta u_x) + \sigma e^{\beta u_x},
\]
where  the parameters $\beta, \sigma,$ and $\lambda$ are different from those in  \eqref{eq:FuxLC1E2}. In fact, the renaming of parameters in a transformed expression to those appearing in the original expression will often be assumed in the sequel. Applying now the equivalence relation \eqref{eq:Fux3} to the latter expression for $F$ reduces it to the form  $F\sim \left[ \left( \delta e^{\alpha u_x} + \lambda \right)/\alpha \right] + \lambda u_x$ for some new arbitrary constants $\delta, \alpha, \lambda.$ Finally inserting the latter expression for $F$ into the determining equations \eqref{eq:detFux} and solving shows that the only existing symmetries for Equation \eqref{eq:Fux:m3E4}  are those for $F$ arbitrary and given by \eqref{eq:FuxArb:sym}.\par

The symmetries corresponding to the third order {\sc ode} \eqref{eq:Fux:m3E5} are found with the same procedure as above for \eqref{eq:Fux:m3E4}. The linear dependence constraints on $F$ takes the form
\begin{subequations}\label{eq:FuxLC2}
\begin{align}
0 &= \alpha u_x + \beta F +\sigma u_x F_{u_x},\label{eq:FuxLC2E1}\\[-2mm]
\intertext{with solution\vspace{-2mm}}
F &= \delta\, u_x^{-\frac{\beta}{ \sigma}} - \frac{\alpha\,  u_x}{\beta + \sigma},  \label{eq:FuxLC2E2}
\end{align}
\end{subequations}

\begin{table}[h]
\caption{ \label{tb:classif1}  \protect { \footnotesize \bf   Nonequivalent symmetry classes for $u_{xy}=F(u_x)\colon$ First five classes.} \footnotesize The $k_j$ are the free parameters of the symmetry group. The functions $\xi= \xi(x), H=H(x)$ and $P=P(y), S=S(y), g=g(y)$ are arbitrary while $\alpha, \beta, \sigma,$ and $\delta$ are free parameters defining the function $F.$ }
%
\begin{tabular}{l l} \hline \\[-1.5mm]
{\scriptsize \bf Numbering }& {\scriptsize \bf Function $F_j$  and Generic symmetry vector $\bv_j,$ for $j=1,\dots, 5.$} \\[1.5 pt]  \hline
$F_1  $  &  $\sim u_x + \alpha,\quad (\alpha =0 \text{ or } \alpha =  1)$     \\[1.5mm]
$\bv_1$&  $\xi \pd_x + (k_1+ g)\pd_y +  \left[g u + \alpha (x g-\xi)+ e^y H(x) + S(y)\right]\pd_u$    \\[3mm]
$F_2$&    $\sim  u_x^\alpha,\quad  \alpha \notin \set{1,2}$    \\[1.5mm]
$\bv_2$&  $(k_3+k_4 x)\pd_x   + (k_1+k_2 y)\pd_y + \left[P+u \left(k_4+\frac{k_2}{1-\alpha}\right) \right] \pd_u $    \\[3mm]
$F_3$&  $\sim  e^{\alpha u_x}\, \sim  e^{ u_x}$    \\[1.5mm]
$\bv_3$&  $(k_2+k_3 x)\pd_x  +   (k_1-k_5 y)\pd_y  +  (k_4+P+k_3 u+k_5 x)\pd_u$    \\[3mm]
$F_4$&  $\sim  u_x^2 + \beta u_x + \beta^2/4\; \sim  (u_x +\alpha)^2,\quad (\alpha =0 \text{ or } 1)$   \\[1.5mm]
$\bv_4$&  \parbox[t]{\textwidth}{$(k_1+k_2 x)\pd_x + [k_3+y (k_4+k_5 y)] \pd_y$ \\[1mm]
 $+ [P+u (k_2-k_4-2 k_5 y)-\frac{1}{2} x (\beta k_4 + 2 k_5+2 \beta k_5 y)]\pd_u$
 } \\[8mm]
$F_5$&    $\sim  u_x^2 + \beta u_x + \beta^2/4 + \delta^2\; \sim u_x(u_x+1) + 1$    \\[1.5mm]
$\bv_5$&  \parbox[t]{\textwidth}{
$(k_1+k_2 x)\pd_x +  \left[\frac{2 \delta k_5-k_4 \cos(2 \delta y)+k_3 \sin(2 \delta y)}{2 \delta}\right]\pd_y $ \\[1mm]
$+ \bigg[ P+k_2 u-\frac{1}{2} (2 k_3
u+\beta k_3 x+2 \delta k_4 x) \cos(2 \delta y) $\\[1mm]
$-\frac{1}{2} (2 k_4 u-2 \delta k_3 x+\beta k_4 x) \sin(2 \delta y)\bigg]\pd_u $
 }   \\[1.5mm]   \hline
\end{tabular}
\end{table}

Thanks to the equivalence relation \eqref{eq:Fux3}, one has $F \sim u_x+ u_x^\beta,$ and substituting the latter expression into \eqref{eq:detFux} gives rise to two relevant cases to consider, namely the case $\beta= 2$ and $\beta \notin \set{1,2}.$ More precisely, the generic symmetry vector $\bv$ corresponding to $F$  is given in this case as follows.
\begin{subequations}\label{eq:Fux:m3SE4}
\begin{alignat}{2}
 F & \sim u_x^\beta + u_x, (\beta \notin \set{1,2}),  &\qquad             \bv &=  V_{3, 6}      \label{eq:Fux:m3SE4a} \\
 F & \sim u_x^2 +u_x,                                  &\qquad             \bv &=  V_{3, 7}, \label{eq:Fux:m3SE4b}
\end{alignat}
\end{subequations}
where
\begin{subequations}\label{eq:Fuxsym3E4}
 \begin{align}
V_{3,6} &=  (k_1 + k_2 x)\pd_x + \left[k_4 + \frac{k_3 e^{y(1-\beta)}}{1-\beta}\right]\pd_y + \left[k_2 +\frac{k_3 e^{y(1-\beta)}}{1-\beta} u \right]\pd_u \label{eq:Fuxsym3E4a}\\
\begin{split}
V_{3,7} &= (k_1 + k_2 x)\pd_x + (e^y k_3 - e^{-y}k_4 + k_5)\pd_y \\
          &\;  + \left(P(y)- k_3\, x e^{y}  + (k_2 - e^y k_3 - e^{-y} k_4)u\right)\pd_u.  \label{eq:Fuxsym3E4b}
\end{split}
 \end{align}
\end{subequations}

\begin{table}[h]
\caption{ \label{tb:classif2}  \protect { \footnotesize \bf   Nonequivalent symmetry classes for $u_{xy}=F(u_x)\colon$ Last five classes.} \footnotesize The $k_j$ are the free parameters of the symmetry group. The functions $\xi= \xi(x), H=H(x)$ and $P=P(y), S=S(y), g=g(y)$ are arbitrary while $\alpha, \beta, \sigma,$ and $\delta$ are free parameters defining the function $F.$ }
%
\begin{tabular}{l l}\hline \\
{ \scriptsize \bf Numbering }& {\scriptsize \bf Function ($F_j$)  and Generic symmetry vector ($\bv_j$)} \\[1.5 pt]  \hline
$F_6$& $\sim  u_x^2 + \beta u_x + \beta^2/4 - \delta^2 \; \sim u_x(u_x+1) $     \\[1.5mm]
$\bv_6$&    \parbox[t]{\textwidth}{
$(k_1+k_2 x)\pd_x +\left[ k_5+\frac{e^{2 \delta y} k_3-e^{-2 \delta y} k_4 }{ 2 \delta} \right]\pd_y$\\[1mm]
$+ \big[P+k_2 u+\frac{1}{2} e^{-2 \delta y} k_4 (-2 u+2 \delta x-\beta x)$\\[1mm]
$ -(1/2)e^{2 \delta y} k_3 (2 u+(2 \delta+\beta) x)\big]\pd_u$
}\\[17mm]
$F_7$&  $\sim  e^{\alpha u_x}+1\; \sim  e^{ u_x}+1 $    \\[1.5mm]
$\bv_7$&  $(k_1 + k_2 x) \pd_x + \left( k_4 - \frac{k_3 e^{- \alpha y}}{\alpha} \right)\pd_y + \left(P(y) + k_2 u - \frac{k_3 x e^{- \alpha y}}{\alpha} \right)\pd_u$    \\[3mm]
$F_8$&  $\sim  (\alpha + u_x)^\beta\; \sim  (1 + u_x)^\beta, \quad  \beta \notin \set{1,2}, \alpha \neq 0 $    \\[1.5mm]
$\bv_8$&   \parbox[t]{\textwidth}{
$(k_2 + k_3\, x)\pd_x +  \left[  k_1 - (\beta -1) \frac{k_5\, y}{\alpha} \right]\pd_y$ \\[1mm]
    $\quad+ \big[ k_4 + k_5\, x + P(y) + \left(k_3 + \frac{k_5}{\alpha}\right)u   \big]\pd_u$
}\\[12mm]
$F_9$&  $\sim  u_x^\beta + u_x,\quad  \beta \notin \set{1,2}$    \\[1.5mm]
$\bv_9$&  $(k_1 + k_2 x)\pd_x + \left[k_4 + \frac{k_3 e^{y(1-\beta)}}{1-\beta}\right]\pd_y + \left[k_2 +\frac{k_3 e^{y(1-\beta)}}{1-\beta} u \right]\pd_u $    \\[3mm]
$F_{10}$&  $F$ is arbitrary.\\[1.5mm]
$\bv_{10}$&  $(k_1 x + k_2)\pd_x + k_3\pd_y + ( k_1 u + P)\pd_u$    \\[1.5mm]   \hline
\end{tabular}
\end{table}

\subsubsection{Four indeterminate are linearly dependent.}
Assuming now that  $m=4$ of the indeterminate in \eqref{eq:detFux} are linearly dependent gives rise to a maximum of four possible conditions on the function $F,$  only two of which yield meaningful solutions and are given as follows.
\begin{subequations}\label{eq:Fux:m4E}
  \begin{align}
    0 &= -F_\nu^2 F_{\nu\nu} + 2 F F_{\nu\nu}^2 - F F_\nu F_{\nu\nu\nu}  \label{eq:Fux:m4E1}\\
    \begin{split}
    0 &= -\nu F_{\nu\nu}^2 F_{\nu\nu\nu} - 3 F F_{\nu\nu\nu}^2 + 2 F F_{\nu\nu}F_{\nu\nu\nu\nu}\\
      &\; + F_{\nu} (2 \nu F_{\nu\nu\nu}^2 + F_{\nu\nu} (F_{\nu\nu\nu}- \nu F_{\nu\nu\nu\nu})) \label{eq:Fux:m4E2}
    \end{split}
  \end{align}
\end{subequations}
Equation \eqref{eq:Fux:m4E1} is just Equation \eqref{eq:GivFu2} whose general solution in terms of $u_x$ is $F\sim (\alpha + u_x)^\beta + \delta,$ and the symmetries associated with this solution for $\delta=0$ are given in \eqref{eq:Fux:m3Sym3}, while the symmetries associated with the case $\alpha =0$ are given in \eqref{eq:Fux:m3E2Sym}. We thus have to assume that $\alpha \neq 0$ and $\delta \neq 0.$ In particular, we may assume that $\delta=1.$  Then while there is no new symmetry associated with the solution $F\sim (\alpha + u_x)^\beta + \delta$ when $\beta \notin \set{1,2},$ there is however one new symmetry $\bv = V_{4,1}$ associated with $F\sim (\alpha+ u_x)^2 + 1,$ $\alpha \neq 0,$ and given by
\begin{equation} \label{eq:Fux:m4Sym}
\begin{split}
 V_{4,1} &=  (k_1 + k_2 x)\pd_x + \left[ k_5 - \frac{1}{2} k_4 \cos(2 y) + k_3 \cos(y) \sin(y)\right] \pd_y \\
           &\quad+ \bigg[P(y) + k_3 x \sin(2 y) - \alpha k_4 x \sin(2y) + k_2 u - k_4 \sin(2y) u\\
           &\quad - \cos(2y)(\alpha k_3 x+ k_4 x + k_3 u)\bigg]\pd_u.
 \end{split}
\end{equation}

On the other hand, Equation \eqref{eq:Fux:m4E2} has general solution
\begin{subequations}\label{eq:m4E2Sol}
\begin{align}
  F &=  \lambda (\tau + u_x)^\beta+ \frac{\alpha\, (\tau + \beta u_x)}{\beta-1} \label{eq:m4E2Sol1}   \\
\intertext{ and by  \eqref{eq:Fux3} one has}
  F &\sim  \alpha\, + u_x + (1+ u_x)^\beta.   \label{eq:m4E2Sol2}
\end{align}
\end{subequations}
Inserting  the expression for $F$ from \eqref{eq:m4E2Sol2}  into the determining equation \eqref{eq:detFux} shows that the associated generic symmetry vector $\bv$ is determined as follows for relevant values of $F.$
\begin{subequations}\label{eq:Fux:m4eq2}
\begin{alignat}{2}
 F & \sim \alpha + u_x + (1+ u_x)^2, (\alpha \neq 5/4),                    &\qquad             \bv &=  V_{4, 2}, \label{eq:Fux:m4eq2b} \\
 F & \sim \frac{5}{4} + u_x + (1+ u_x)^2,                                  &\qquad             \bv &=  V_{4,3}, \label{eq:Fux:m4eq2c}
\end{alignat}
\end{subequations}
where
\begin{subequations}\label{eq:Fuxsym4E2}
 \begin{align}
 \begin{split}\label{eq:Fuxsym4E2a}
V_{4,2} &=  (k_1 + k_2 x)\pd_x  + \left[k_5 + \frac{k_3 e^{y \sqrt{5-4 \alpha}}   }{  \sqrt{5-4 \alpha}}  -  \frac{k_4 e^{-y \sqrt{5-4 \alpha}}   }{  \sqrt{5-4 \alpha}}  \right]  \pd_y  \\
&\quad +\quad  \bigg[ \frac{1}{2}   e^{-y \sqrt{5 -4 \alpha}} \big(2 e^{y \sqrt{5 -4 \alpha}} S(y) + k_4 (x(-3 + \sqrt{5 -4 \alpha})  -2 u) \\
&\quad + 2 k_2 e^{y \sqrt{5 -4 \alpha}} u -  k_3 e^{2 y \sqrt{5 -4 \alpha}}\left( x(3+ \sqrt{5 -4 \alpha} + 2 u)\right)    \big) \bigg]\pd_u
\end{split}\\
\begin{split}
V_{4,3} &= (k_1 + k_2 x)\pd_x + (k_3 + y(k_4+ k_5 y))\pd_y \\
          &\;  + \left(S(y)-\frac{1}{2} x (3 k_4 + 2 k_5 + 6 k_5 y ) + (k_2 - k_4 - 2 k_5 y)u     \right)\pd_u.  \label{eq:Fuxsym4E2b}
\end{split}
\end{align}
\end{subequations}

\subsubsection{Five indeterminate are linearly dependent.}
We now consider the last case where in \eqref{eq:detFux} there are five linearly dependent indeterminate in $u_x$. The corresponding single equation is given by
  \begin{equation}\label{eq:Fu:m5E}
    0 = -F_{\nu\nu\nu}^2 F_{\nu\nu\nu\nu} + 2 F_{\nu\nu} F_{\nu\nu\nu\nu}^2 - F_{\nu\nu}F_{\nu\nu\nu}F_{\nu\nu\nu\nu\nu}.
  \end{equation}
Setting  $F= \int \int G(\nu) d\nu d\nu,$ and then $G= e^{w(\nu)}$ reduces \eqref{eq:Fu:m5E} to the simpler equation $2w_{\nu\nu}^2 - w_{\nu}w_{\nu\nu\nu}$, whose solution is $w= \beta \ln(\alpha + \nu) + \delta.$ However, none of corresponding solutions of the original equation \eqref{eq:Fu:m5E} yields a new symmetry.\par

Before concluding our analysis for the case $F=F(u_x),$ it should be noted that functions $F$ found in \eqref{eq:Fux:m2Sym1},   \eqref{eq:Fux:m2Sym2}, \eqref{eq:Fux:m3Sym3a}, \eqref{eq:Fux:m3SE4b}, \eqref{eq:Fux:m4Sym}, \eqref{eq:Fux:m4eq2} all satisfy $F \sim u_x^2 + \beta u_x + \alpha$ for some constant parameters $\alpha$ and $\beta,$ and  are therefore equivalent to one of the three functions classified earlier in \eqref{eq:Fux:m3Sym1}.

\begin{thm}\label{th:classFux}

For the hyperbolic equation $u_{xy}= F(u_x),$  the list of all possible values of $F=F(u_x)$ yielding, under point transformations,  nonequivalent symmetry algebras   with generic vectors $\bv$ are given by Table \ref{tb:classif1} and Table \ref{tb:classif2}, and consists of a total of ten  nonequivalent symmetry classes.\par
\end{thm}

\begin{proof}
The fact that Table \ref{tb:classif1} and Table \ref{tb:classif2} list all possible distinct  values of $F$ and the corresponding symmetry algebra  follows from the discussion carried out in Section \ref{subsecFux}. The simplification of the canonical form of each function $F$ listed in the table naturally follows from \eqref{eq:Fux3}. It remains however to ascertain that all ten functions in these lists are pairwise nonequivalent. This is clear for all functions and separately for the group  of functions $F_4, F_5$ and $F_6,$ in view of the fact that by \eqref{eq:Fux3}, the equivalence transformations \eqref{eq:Fux2} acts of $F$ and its argument by mere scalings. Moreover, the functions $F_4, F_5,$ and $F_6$ are pairwise nonequivalent as this is clearly the case for their canonical forms, namely for the functions $(u_x +1)^2,$ $u_x (u_x +1),$ and $(u_x+1)^2 +1.$
\end{proof}

\section{Conservation laws} \label{s:clws}

As discussed in Section \ref{s:intro}, the many physical applications associated with several equations of the form \eqref{eq:M} motivates the study in this section of conservation laws of this class of equations. Let us recall that a conservation law of \eqref{eq:M}, is a divergence expression
\begin{equation}\label{eq:consLw1}
\diverg \Theta \equiv D_x \Phi + D_y \Psi=0, \quad \Theta = \Theta[u]= (\Phi, \Psi)= (\Phi[u], \Psi[u]),
\end{equation}
that vanishes on the solution space of \eqref{eq:M}, where $D_x$ and $D_y$ are the usual total differential operators acting on the jet space of the underlying space $X\times U \subset \R^3$ of independent variables $(x,y) \in X$ and  dependent variable $u \in U.$ Moreover, the notation $f=f[u]$ denotes a differential function of $u,$ that is, a function depending on $x,y, u$ and the derivatives of $u$ up to an arbitrary but fixed order. The vector $\Theta= (\Phi, \Psi)$ is then called a conserved current, and may also be called a flux vector, as we shall do.\par

Conservation laws are also determined up to an equivalence class, whereby two conservation laws are termed equivalent if they differ by a trivial conservation law, that is by one which vanishes for all smooth functions $u=u(x,y)$ and not only on the solution space of \eqref{eq:M}. Such trivial conservation laws are said to be of the second kind as opposed to trivial conservation laws of the first kind  in which the function $\Theta$  in \eqref{eq:consLw1} vanishes itself on all solutions of the equation  \cite{olv93}. Trivial conservation laws of the first kind will be eliminated by the choice of derivatives appearing as arguments in $\Theta.$ \par

To each conservation law \eqref{eq:consLw1}, there corresponds a multiplier $Q=Q[u]$ given by
\begin{equation}\label{eq:charact}
Q \cdot \Delta =  D_x \Phi[u] + D_y \Psi[u]=0,\quad \Delta \equiv  u_{xy} - F(u u_x).
\end{equation}
and two multipliers are equivalent if they differ by a trivial multiplier, that is, by one which vanishes itself on the solution space of the equation. For normal and totally nondegenerate equations of the form \eqref{eq:M}, there is a one-to-one correspondence between equivalent classes of multipliers  and equivalent classes of conservation laws.\par
For Lagrangian equations, every characteristic of a variational symmetry is a multiplier and can be uniquely identified with  a conservation law. However, Lagrangian equations constitute a very restricted type of equations and for more general equations for which the concept of variational symmetry in particular does not  apply, multipliers are found as generalized integrating factors. That is, as factors $Q=Q[u]$ for which the product $Q\cdot \Delta$ is a total divergence expression, and hence a null Lagrangian.\par

Our search for conservation laws will be focussed on those of low order, as physically relevant properties such as energy, mass-energy, and momentum conservation are always associated with these low-order conservation laws. In mathematical terms, these are conservation laws for which every derivative of a dependent variable in the equation can be obtained by differentiating with respect to some independent variable a similar derivative in the expression of the multiplier.\par

Conservation laws of \eqref{eq:M} with $F=F(u)$ are given for arbitrary values of $F$ in \cite{polyanin1}, although \cite{polyanin1} does not make any reference to the corresponding multiplier. It is indeed possible to find a specific flux vector  $\Theta$ directly, albeit generally more intricate, by solving  \eqref{eq:consLw1} on the solution surface of the equation.  We complement the result of \cite{polyanin1} by finding the multiplier $Q$ which is to be sought in the form $Q= Q(x,y, u_x)$ or $Q(x,y, u_y).$ In view of obvious symmetry considerations associated with the structure of \eqref{eq:Fu}, it is enough to search for $Q= Q(x,y, u_x).$ The determining equations for the function $Q$ are
\begin{align*}
  Q_u &=0     \\
  Q_{y u_x} + F Q_{u_x u_x}&=0         \\
F_u Q_{u_x} u_x -Q F_u - F Q_u + F u_x Q_{u u_x} + Q_{x y} + F Q_{x u_x} + u_x Q_{y u}&=0,
\end{align*}
with solution $Q= \alpha u_x,$ for a certain constant parameter $\alpha,$ which we may assume to be equal to $1$ thanks to the linearity of the Euler Operator. Indeed, $Q$ is determined by the condition $E_u(Q\cdot \Delta)=0,$ where $E_u$ is the Euler Operator with respect to $u.$ The corresponding determining equations for the flux  vector $\Theta=(\Phi, \Psi)$ takes the form
\begin{align*}
0&={\Phi}_{u_{xx}},\qquad{\Phi}_{u_{yy}} =0,\qquad{\Psi}_{u_{xx}}=0,\qquad {\Psi}_{u_{yy}}=0    \\
0&=u_x -{\Phi}_{u_y} - {\Psi}_{u_x}          \\
0&=u_x (-F - {\Phi}_u) - u_{xx}{\Phi}_{u_x} -{\Phi}_x - u_y{\Psi}_u - u_{yy}{\Psi}_{u_y} -{\Psi}_y .
\end{align*}
The latter determining equation has solution
\[
\Theta= \left(-F +\int R_y \, du+ u_y R+S,\; \frac{u_x^2}{2}-\int S_x \, dy -\int R_x
\, du+P-u_x R   \right)
\]
where $R=R(x,y,u), S=S(x,y), P=P(x)$ are arbitrary functions. It turns out that this flux vector contains as a term a trivial flux , that is,  one corresponding to a trivial conservation law, namely
the term
\[
 \left( \int R_y \, du+ u_y R+S,\; -\int S_x \, dy -\int R_x
\, du+P-u_x R   \right).
\]
Using the symmetry argument related to the structure of \eqref{eq:Fu} together with the above results show that the multipliers $Q$ and corresponding conservation laws of \eqref{eq:Fu} are given for $F$ arbitrary by
\begin{subequations} \label{eq:ClwFu}
\begin{alignat}{2}
Q &= u_x,\qquad&   0&= D_x \left[ -F(u)   \right] +   D_y ( u_x^2/2)  \label{eq:ClwFu1}\\
Q &= u_y,\qquad&   0&= D_x ( u_y^2 /2) + D_y \left[ -F(u)  \right].    \label{eq:ClwFu2}
\end{alignat}
\end{subequations}
We shall make use of the same procedure used thus far in this section to find the conservation laws of other relevant equations of the form \eqref{eq:M}. In this way, the next equation from the class \eqref{eq:M} with $F=F(u)$ which we consider is the Liouville equation
\begin{equation}\label{eq:liouv}
  u_{xy} = e^u.
\end{equation}
Although the conservation laws of  \eqref{eq:liouv}, which is a well studied equation \cite{winterLv1, winterLv2, dubrovin} are likely to have been computed, we find it interesting to gather all such results in a single short paper.  After calculations, the multipliers $Q$ and corresponding conservations laws for \eqref{eq:liouv} are given as follows.
\begin{subequations}\label{eq:clwLiouv}
\begin{alignat}{2}
  Q &= p_x+ p\, u_x,\quad &     0 &= D_x\left[ -e^u p + p_x u_y  \right]         + D_y \left[ -u p_{xx}+ (1/2) p\, u_x^2 \right]  \\
  Q &=q_y+ q\, u_y,\quad &     0 &= D_x \left[  -u q_{yy}+ (1/2) q\, u_y^2 \right] + D_y\left[-e^u q + q_y u_x   \right],
\end{alignat}
\end{subequations}
where $p=p(x)$ and $q=q(y)$ are arbitrary functions. This shows in particular that \eqref{eq:liouv} has infinitely many conservation laws of low order. It has in fact been proved, that \eqref{eq:clwLiouv} is symmetry integrable, that is, it has an infinite series of generalized symmetries of arbitrary orders. Moreover \cite{sokolov}, all symmetry integrable equations of the form \eqref{eq:Fu} are given up to scalings and shifts by functions $F=F(u)$ satisfying
\begin{equation}\label{eq:FuExpV2}
F(u)= e^u,\qquad F(u)= e^u + e^{-u},\qquad  \text{ or }\quad  F(u) =  e^u + e^{-2u}.
\end{equation}
We now let $F(u)= u^2 + 1$ in \eqref{eq:Fu}, and the resulting equation takes the form.
\begin{equation}\label{eq:FuSqP1}
u_{xy} = u^2 +1.
\end{equation}
One sees that  low-order multipliers are linear combinations of $u_x$ and $u_y.$ Letting the multiplier $Q= u_x,$ the corresponding flux vector has expression
\begin{equation}\label{eq:FuSqP1_Flux1}
\Theta = \left(B-u-\frac{u^3}{3}+A u_y+  \int A_y  \, du,  C-A u_x+\frac{u_x^2}{2} -\int {B}_x  \, dy -\int A_x  \, du \right),
\end{equation}
where $A= A(x,y,u),$ $B= B(x,y),$ and $C=C(x)$ are arbitrary functions of their arguments.  The pure flux vector associated with the multiplier $Q= u_x,$ that is, the one deprived of any trivial flux vector  is however given  by
\[
\Theta = \left( -u-\frac{u^3}{3},  \frac{u_x^2}{2} \right).
\]
In other words, the aggregate of arbitrary functions in \eqref{eq:FuSqP1_Flux1} only yields a trivial flux vector.  This results and the symmetrical structure of \eqref{eq:Fu} with respect to the variable $x$ and $y$ shows that all multipliers and corresponding conservation laws of \eqref{eq:FuSqP1} are given by
\begin{subequations}\label{eq:clwFuSqP1}
\begin{alignat}{2}
  Q &=u_x,\quad &     0 &= D_x\left( -u- u^3/3 \right) + D_y \left( u_x^2/2  \right)  \\
  Q &=u_y,\quad &     0 &= D_x \left(  u_y^2/2   \right) + D_y\left(  -u- u^3/3  \right).
\end{alignat}
\end{subequations}
In particular, \eqref{eq:FuSqP1} has only two conservation laws of low order.\par

For the counterpart \eqref{eq:Fux} of \eqref{eq:Fu} for which $F= F(u_x)$ so that the equation reduces to $u_{xy} = F(u_x),$  no general result of a similar nature concerning the symmetry integrability of this equation seems to be available. To begin with,  for arbitrary values of $F,$ \eqref{eq:Fux} has only  trivial multipliers and hence trivial conservation laws. Nevertheless, it turns out that for the particular case $F(u_x)= u_x^2,$ the corresponding equation
\begin{equation}\label{eq:Fux^2}
  u_{xy}= u_x^2
\end{equation}
also has an infinity of conservation laws of  low order.

Indeed, using our usual procedure described above, the general expression for the multiplier $Q$ of low-order conservation laws of \eqref{eq:Fux^2} takes the form
\begin{subequations}\label{eq:multUx^2}
  \begin{align}
Q &= T(z, y) +  \frac{x}{z^2} V\left(x, -\frac{1+ z y}{z}\right),\quad z= u_x,  \label{eq:multUx^2v1}\\
\intertext{where $V$ is an arbitrary function of its arguments, and $T= T(z,y)$ is a solution of the second order linear hyperbolic equation}
0&= 2 T + 4 z T_z + z^2 T_{z,z}  + T_{zy}.    \label{eq:multUx^2v2}
  \end{align}
\end{subequations}
Although a general solution of \eqref{eq:multUx^2v2} is not available, we can find its particular solutions, and this will be enough to verify this claim. A practical and relatively simple way to obtain these particular solutions is to use the similarity reduction method yielding group-invariant solutions \cite{olv93, bluman, ovsyannik, stephani}.  In addition to the solution symmetry $\zeta \pd_T$, where $\zeta= \zeta(y,z)$ is any given solution of \eqref{eq:multUx^2v2}, the other symmetries of this equation are given by
\begin{subequations}\label{eq:symMultV1}
\begin{alignat}{2}
\bv_1 &=  T\pd_T,\qquad                                  &     \bv_2&= -z \pd_y + y \pd_z,  \\
\bv_3 &= (-2y z-1)\pd_y + y^2 \pd_z + 2 T y \pd_T,\quad  &     \bv_4&=  \pd_z.
\end{alignat}
\end{subequations}
Group-invariant solutions $T_1, T_2, T3$ associated with the above symmetries $\bv_1, \bv_2, \bv_3$ of \eqref{eq:multUx^2v2} and the corresponding multipliers $Q_1, Q_2,$ and $Q_3$ of \eqref{eq:Fux^2} are given by
\begin{subequations}\label{eq:Q1to3}
\begin{alignat}{2}
T_1 &= \frac{\beta_1}{ 4 z^2} -\frac{\alpha_1}{2 z} ,\qquad&    Q_1 &= T_1+ \frac{x}{u_x^2} V \\
T_2 &= \frac{\alpha_2}{(1+ y z)^2}  + \frac{\beta_2(y z + \ln (y z))}{(1+ y z)^2},\qquad&    Q_2 &=  T_2 +  \frac{x}{u_x^2} V \\
T_3 &= \frac{\alpha_3-y - y^2 u_x}{\beta_3(1+ y u_x)^2} ,\qquad&    Q_3&= T_3 +  \frac{x}{u_x^2} V,
\end{alignat}
\end{subequations}
where the $\alpha_j$ and $\beta_j$ are arbitrary constants for $j=1,2,3,$ and the arbitrary function $V$ is given by \eqref{eq:multUx^2v1}. Finally, the flux vectors $\Theta_j=(\Phi_j, \Psi_j)$ corresponding to the multipliers $Q_j$ are given by
\begin{subequations}\label{eq:fluxUx^2}
\begin{align}
  \Theta_1    &=  \left(\frac{\alpha_1 u}{2} - \frac{\beta_1 x}{4},\quad  - \frac{1}{4} \left(\frac{\beta_1}{u_x}  + 2 \alpha_1 \ln(u_x) \right)  + \int \frac{x}{u_x^2} V d u_x\right)  \\
  \Theta_2    &= \left[ \frac{-\beta_2 u}{y},\;  \frac{\alpha_2 - 2 \beta_2}{y} +  \int \frac{\alpha_2 + \beta_2 u_x y + \beta_2 \ln (u_x y)}{(1+ u_x y)^2} +  \frac{x}{u_x^2} V d u_x  \right]  \\
  \Theta_3    &= \left[  \frac{u}{\beta_3},\;  \frac{\beta_3^{-1}}{ y}
    \left( \alpha_3 - y \ln(y) + y \int \left( \frac{a-y  - u_x y^2}{(1+ u_x y)^2}+ \frac{x}{u_x^2} V \right) d u_x   \right)\right].
\end{align}
\end{subequations}
Clearly, the nontrivial conservation law for each of the pure flux vectors $\Theta_j= (\Phi_j, \Psi_j)$ above are given by $D_x \Phi_j + D_y \Psi_j=0,$ and such low-order conservation laws are infinitely many as they dependent of the arbitrary function $V= V(x, p),$ with $p= - (1+ u_x y)/ u_x.$


\section{Concluding Remarks}
\label{s:conclusion}
The Lie group classification of differential equations has proved to be a very challenging exercise, even for the simplest types of differential equations that {\sc ode}s represent, and we have made a breakthrough in this paper by finding a simple and systematic way of solving this type of problems, based on the determination of all possible cases of linearly dependent indeterminate in the determining equations.  This has allowed us to give a complete classification of the family \eqref{eq:M} of equations in a relatively brief manner. Although the method has only been discussed in application to the particular case of Equation  \eqref{eq:M}, it seems  however quite clear that it can be extended in a straightforward way to at least any scalar differential equation, and in fact to any system of {\sc ode}s or {\sc pde}s. In the case of systems of equations, it might however be more prudent to make a more confident conclusion only after an adequate number of concrete applications of the method.\par

Calculations done in Section \ref{s:clws}
show in particular that  Equation \eqref{eq:FuSqP1} has only a finite number of low-order conservation laws, while \eqref{eq:liouv} and \eqref{eq:Fux^2} have infinitely many of them. However, while  \eqref{eq:liouv} is known to be up to two other related cases listed in \eqref{eq:FuExpV2} the only equations of the form  \eqref{eq:Fu} that are symmetry integrable, nothing is known about the symmetry integrability of  \eqref{eq:Fux^2} and in particular about that for the more general case \eqref{eq:Fux}. The symmetry integrability property of differential equation remains a very active and challenging domain of research with very limited results, and it would be interesting to find out if \eqref{eq:Fux^2} is symmetry integrable, and if so, whether such property also relates to the symmetry integrability of the whole family of equations \eqref{eq:Fux}.

\section*{Declarations of interest:} None.
\label{s:declare}

\section*{Acknowledgement}
\label{s:acknowledge}
Funding: This work was supported by the NRF Incentive Funding for Rated Researchers grant [Grant Number 97822];
the University of Venda [Grant Number I538];

\bibliographystyle{model1-num-names}

\end{document}